\documentclass[a4paper,11pt]{article}
\usepackage{multirow}
\usepackage{pstricks}
\usepackage{pst-node}
\usepackage{pst-tree}
\usepackage{pst-eps}
\usepackage{xspace}

\usepackage{amsfonts}

\usepackage{latexsym}
\usepackage{amsthm}
\usepackage{fancyheadings}
\usepackage{amsmath}
\usepackage{graphicx}

\newtheorem{Proposition}{Proposition}
\theoremstyle{definition}

\title{A trace partitioned Gray code for 
$q$-ary generalized Fibonacci strings} 
\author{A. Bernini, S. Bilotta, R. Pinzani \\ 
        {\it\small Dipartimento di Matematica e Informatica ``Ulisse Dini"}\\
        {\it\small Universit\`a degli Studi di Firenze, Viale G.B. Morgagni 65, 50134 Firenze, Italy}\\ 
        {\tt\small  e-mail: \{antonio.bernini\}\{stefano.bilotta\}\{renzo.pinzani@unifi\}@unifi.it}\\ \\
	V. Vajnovszki\\
	{\it\small LE2I, Universit\'e de Bourgogne, BP 47 870, 21078 Dijon Cedex, France}\\
	{\tt\small  e-mail: vvajnov@u-bourgogne.fr}
	}

\def\first{{ \mathsf{first}}}
\def\last{{ \mathsf{last}}}
\begin{document}

\maketitle

\begin{abstract} 
\noindent
We provide a trace partitioned Gray code for the set of 
$q$-ary strings avoiding a pattern constituted by $k$ consecutive
equal symbols. 
The definition of this Gray code is based on two different
constructions, according to the parity of $q$.
This result generalizes, and is based on, a Gray code for binary strings avoiding 
$k$ consecutive~0's.
\end{abstract}

{\it Keyword}:
Gray codes, 
pattern avoiding strings, 
generalized Fibonacci numbers

\section{Introduction} 
The famous \emph{$k$-generalized Fibonacci sequence}
$\{f_n^{(k)}\}_{n\geq0}$ is defined by
$$
f_n^{(k)}=f_{n-1}^{(k)}+f_{n-2}^{(k)}+\cdots+f_{n-k}^{(k)},\mbox{ for}\ n\geq k,
$$
with initial conditions \cite{K}:
$
f_i^{(k)}=0$ for $0\leq i\leq k-2$, and $f_{k-1}^{(k)}=1.
$
This sequence is related to the 
enumeration of binary strings avoiding $k$ consecutive 1's \cite{K},
called \emph{$k$-generalized Fibonacci strings}.

If an alphabet of cardinality $q\geq 2$ is used, then the enumeration of the strings of
length $n$ and avoiding a pattern constituted by $k$ consecutive occurrences of a same symbol
is given by

$$
f_{n,q}^{(k)}=(q-1)\left(f_{n-1,q}^{(k)}+f_{n-2,q}^{(k)}+\ldots+f_{n-k,q}^{(k)}\right),
\mbox{ for}\ n\geq k,
$$
with
$f_{i,q}^{(k)}=q^i$, for $0\leq i\leq k-1$;
and in particular $f_{n+k}^{(k)}=f_{n,2}^{(k)}$.
This sequence is a particular case of the \emph{weighted $k$-generalized Fibonacci sequence}
(studied and used in \cite{L,Sh}) where all the weights are equal to $q-1$.
Similarly to the binary case, by a
\emph{$q$-ary $k$-generalized Fibonacci string} we mean 
a string over a $q\geq 2$ letter alphabet $A$ and avoiding $k$ consecutive occurrences of a 
given symbol of $A$. 
For example, for $k=3$ and $A=\{0,1,2,3\}$,
$111$ or $222$ is avoided, but not both.

Once a class of objects is defined, often it could be of interest to list or generate them according 
to a particular criterion.
A special way to do this is to generate the objects so that any two consecutive ones
differ as little as possible, i.e., in a \emph{Gray code} manner \cite{G}. Such a particular 
code has already been proposed for (binary) $k$-generalized Fibonacci strings \cite{V1}.

In the present work we provide a Gray code for the
$q$-ary $k$-generalized Fibonacci strings, $q\geq 2$, which extends the approach in \cite{V1}. 
Our method gives a \emph{trace partitioned} code where strings with the same trace are 
contiguous (for more details see Section \ref{Vajno}).
Analogously to the kind of the forbidden pattern considered in \cite{V1} 
(which is a sequence of $k$ consecutive $1$'s), 
here we consider the avoidance of $k$ consecutive occurrences of a give
symbol, not necessarily $1$, of a $q$-ary alphabet.

\medskip
If $L$ is a set of strings over an alphabet $A$, then by $\mathcal L$ we denote the ordered list
where the strings of $L$ are listed following a certain criterion. If the Hamming distance
\cite{H} between two successive elements of $\mathcal L$ is bounded by a constant,
than $\mathcal L$ is called Gray code list. The notations we are going to use 
are defined below:
\begin{itemize}
\item[-] if $\alpha$ and $\beta$ are two same length strings,
then $d_H(\alpha,\beta)$ denotes their Hamming distance;
\item[-] $\overline{\mathcal L}$ denotes the list obtained by covering 
 $\mathcal L$ in reverse order;
\item[-] $\first(\mathcal L)$ and $\last(\mathcal L)$ are the first and the last element
        of $\mathcal L$, respectively, and clearly $\first(\mathcal L)=\last(\overline{\mathcal L})$,
	and $\last(\mathcal L)=\first(\overline{\mathcal L})$;
\item[-] if $u$ is a string, then $u\cdot\mathcal L$ is the list 
        obtained by prepending $u$ to each string in $\mathcal L$;
\item[-] if $\mathcal L$ and $\mathcal L'$ are two lists, then $\mathcal L\circ\mathcal L'$
is the concatenation of the two lists, obtained by appending the elements of
$\mathcal L'$ after the elements of $\mathcal L$.
\end{itemize}

\section{Gray codes}\label{Vajno}

The well known Binary Reflected Gray Code (BRGC) \cite{G} can be generalized to
an alphabet of cardinality greater than 2 \cite{E,Wi}. If $A=\{a_0,a_1,\ldots,a_{q-1}\}$,
then the list $\mathcal G_n^q$ of the length $n$ strings over $A$ is given by:

\begin{equation}
\mathcal G_n^q=
\left\{
\begin{array}{cr}
\lambda & \quad\mbox{if}\quad n=0,\\
\\
a_0\cdot\mathcal G_{n-1}^q\circ a_1\cdot\overline{\mathcal G_{n-1}^q}
\circ\cdots\circ a_{q-1}\cdot\mathcal G_{n-1}^{'q}&\quad\mbox{if}\quad n>0,
\end{array}
\right.
\label{G_C_General}
\end{equation}

\noindent where $\mathcal G_{n-1}^{'q}$ is $\overline{\mathcal
G_{n-1}^q}$ if $q$
is even or $\mathcal G_{n-1}^q$ if $q$ is odd, and $\lambda$ is the empty string.
It is proved that $\mathcal G_n^q$ is a Gray code list with Hamming distance 1 and the reader can easily prove the following proposition.

\begin{Proposition}\label{first_last}

If $q$ is odd, then $\last(\mathcal G_n^q)=a_{q-1}^n$ and
$\first(\mathcal G_n^q)=a_0^n$. In particular,
$$\last(\mathcal G_n^q)=\last(\mathcal G_{n-1}^q)\cdot a_{q-1}=a_{q-1}\cdot\last(\mathcal G_{n-1}^q),$$
and
$$\first(\mathcal G_n^q)=\first(\mathcal G_{n-1}^q)\cdot a_0=a_0\cdot \first(\mathcal G_{n-1}^q).$$\\
If $q$ is even, then $\last(\mathcal G_n^q)=a_{q-1}a_0^{n-1}$ and
$\first(\mathcal G_n^q)=a_0^n$. In particular, for $n>1$, 
$$\last(\mathcal
G_n^q)=a_{q-1}a_0^{n-1}=\last(\mathcal G_{n-1}^q)\cdot a_0,$$
and
$$\first(\mathcal G_n^q)=a_0^n=\first(\mathcal G_{n-1}^q)\cdot a_0=a_0\cdot 
\first(\mathcal G_{n-1}^q).$$
\end{Proposition}

We will give a Gray code for the set of $q$-ary $k$-generalized 
Fibonacci strings of length $n$, where
the Hamming distance between two consecutive strings is~$1$.
The definition of our Gray code depends on the parity of $q$, and 
we will refer simply to Gray code for a list where two successive elements 
have Hamming distance equal to $1$, and without any loss of generality 
we will consider the alphabet $A=\{0,1,\ldots,q-1\}$. For our topics we are going to 
consider the avoidance of the pattern $0^k$, for a given $k\geq 2$,
but all of our constructions can be easily translated to the pattern $i^k$, for any 
$i\in A$.

Our approach starts from the definition of the Gray code for binary $k$-generalized 
Fibonacci strings given in \cite{V1}, then using a bit replacing technique we extend 
the binary alphabet to $A$, leading to a trace partitioned Gray code
(where strings with the same trace are consecutive). 
The \emph{trace} of a $q$-ary string is a binary string obtained
by replacing each symbol different from 0 by~1.

Before going along in our discussion, we define a tool for
manipulating the symbols of a string. 
For $q\geq 3$, we denote by $\mathcal G_t^{q-1}\oplus 1$ the list 
obtained from $\mathcal G_t^{q-1}$ by adding $1$ to each symbol in each string in 
$\mathcal G_t^{q-1}$. Actually, $\mathcal G_t^{q-1}\oplus 1$ is the 
Gray code defined in relation (\ref{G_C_General}) for 
$\{a_0,a_1,\ldots,a_{q-1}\}=\{1,2,\ldots,q\}$.
For example,  $\mathcal G_3^2\oplus 1=(111,112,122,121,221,222,212,211)$.

If $\beta$ is a binary string of length $n$ such that $|\beta|_1=t$ (the number of 1's in
$\beta$), we define the \emph{expansion} of $\beta$, denoted by $\varepsilon(\beta)$, 
as the list of $(q-1)^t$ strings, where the $i$-th string is obtained by replacing
the $t$ 1's of $\beta$ by the $t$ symbols (read from left to right) of the $i$-th string in
$\mathcal G_t^{q-1}\oplus1$.
For example, if $q=3$ and $\beta=01011$ (the trace),  then
with $\mathcal G_3^2\oplus 1$ given above, we have 
$\varepsilon(\beta)=(01011,01012,01022,01021,02021,02022,02012,02011).$
Notice that in particular $\first(\varepsilon(\beta))=\beta$ and all the elements of
$\varepsilon(\beta)$ have the same trace.

We observe that $\varepsilon(\beta)$ is the list obtained from 
$\mathcal G_t^{q-1}$ by adding 1 to each symbol of each string in $\mathcal G_t^{q-1}$,
then inserting
some 0's, each time in the same positions.
Since $\mathcal G_t^{q-1}$
is a Gray code and the insertions of the 0's does not change the
Hamming distance between two successive element of
$\varepsilon(\beta)$ (which is 1), we have the following:

\begin{Proposition}\label{epsilon}
For any $q\geq 3$ and binary string $\beta$, the list $\varepsilon(\beta)$ is a Gray code.
\end{Proposition}

The Gray code we are going to consider as the starting point of
our argument is the one defined in \cite{V1}, where the author
deals with binary strings avoiding $k$ consecutive 1's. Since we
are interested in the avoidance of $k$ consecutive 0's, we recall,
for the sake of clearness, the definition in \cite{V1} adapted
according to our needs which cause some slight differences with
respect to the original definition in \cite{V1}.

Let $\mathcal{F}_n^{(k)}$ be the list defined as:

{\small
\begin{equation}
\label{f_nk}
\mathcal{F}_n^{(k)}=\left\{
\begin{array}{cr}
\mathcal C_n &\ \mathrm{if}\ 0\leq n<k,\\
\\
1\cdot\overline{\mathcal F_{n-1}^{(k)}}\circ 01\cdot\overline{\mathcal F_{n-2}^{(k)}}
\circ 001\cdot\overline{\mathcal F_{n-3}^{(k)}}\circ\ldots\circ 0^{k-1}1\cdot
\overline{\mathcal F_{n-k}^{(k)}}& \ \mathrm{if}\ n\geq k,
\end{array}
\right.
\end{equation}
}
where
$$
\mathcal C_n=
\left\{
\begin{array}{cr}
\lambda&\quad\mbox{if}\quad n=0,\\
\\
1\cdot\overline{\mathcal C_{n-1}}\circ 0\cdot\mathcal C_{n-1}&\quad\mbox{if}\quad n\geq 1,

\end{array}
\right.
$$
with $\lambda$ the empty string.

\medskip

It is proved \cite{V1} that $\mathcal{F}_n^{(k)}$ is a 
Gray code for the set of binary length $n$ strings avoiding 
$k$ consecutive $0$'s, and thus the number of strings in 
$\mathcal{F}_n^{(k)}$ is $f_{n+k}^{(k)}$, the $(n+k)$-th value of the 
$k$-generalized Fibonacci sequence.

Now let $F_{n,q}^{(k)}$ be
the set of length $n$ strings over $A=\{0,1,\ldots,q-1\}$, $q>2$, avoiding $k$ consecutive 0's,
so $|F_{n,q}^{(k)}|=f_{n,q}^{(k)}$ (see Introduction). The aim is the
construction of a Gray code for $F_{n,q}^{(k)}$.
Our definition of such a Gray code is based on 
the expansion $\varepsilon(\alpha_i)$
(or its reverse $\overline{\varepsilon(\alpha_i)}$)
of each element of $\mathcal{F}_n^{(k)}=\left(\alpha_1,\alpha_2,\ldots,\alpha_{f_{n+k}^{(k)}}\right)$,
and then concatenating them opportunely, according to the parity of $q$.

Let us illustrate the construction of the list for the set $F_{n,q}^{(k)}$,
of the form 

\begin{equation}
\label{Gent_Alter}
\varepsilon(\alpha_1)\circ\overline{\varepsilon(\alpha_2)}
\circ\varepsilon(\alpha_3)\circ\overline{\varepsilon(\alpha_4)}
\cdots
\end{equation}
where $\alpha_i$ covers the list $\mathcal F_n^{(k)}$.
As we will see below, this construction yields a Gray code when
$q$ is even, but not necessarily when $q$ is odd.

When  $k=q=n=3$, the list  $\mathcal F_n^{(k)}$ is 
$\mathcal F_3^{(3)}=~(100,101,111,110,010,011,001)$;
$\mathcal G_1^2\oplus 1=(1,2)$, $\mathcal G_2^2\oplus 1=(11,12,22,21)$, and $\mathcal G_3^2\oplus 1$
is given in the example preceding Proposition \ref{epsilon}. 
The expansions of the elements of  
$\mathcal F_3^{(3)}$ are:
$$
\begin{array}{lllll}
\varepsilon(\alpha_1)&=&\varepsilon(100)&=&(100,200)\\
\varepsilon(\alpha_2)&=&\varepsilon(101)&=&(101,102,202,201)\\
\varepsilon(\alpha_3)&=&\varepsilon(111)&=&(111,112,122,121,221,222,212,211)\\
\varepsilon(\alpha_4)&=&\varepsilon(110)&=&(110,120,220,210)\\
\varepsilon(\alpha_5)&=&\varepsilon(010)&=&(010,020)\\
\varepsilon(\alpha_6)&=&\varepsilon(011)&=&(011,012,022,021)\\
\varepsilon(\alpha_7)&=&\varepsilon(001)&=&(001,002),\\
\end{array}
$$
and the list of the form (\ref{Gent_Alter}) for $F_{3,3}^{(3)}$ is
$$\varepsilon(\alpha_1)\circ\overline{\varepsilon(\alpha_2)}
\circ\varepsilon(\alpha_3)\circ\overline{\varepsilon(\alpha_4)}
\circ\varepsilon(\alpha_5)\circ\overline{\varepsilon(\alpha_6)}
 \circ\varepsilon(\alpha_7)$$ 
which is a Gray code, as it easily can be checked.
However, this is not true in general.
For example, if $n=4$, $k=q=3$, concatenating $\epsilon(\alpha_i)$ and 
$\overline{\varepsilon(\alpha_{i+1})}$, alternatively as in 
(\ref{Gent_Alter}) does not yield a Gray code.
Indeed, 
\begin{itemize}
\item $\alpha_7=1100$ and $\alpha_8=0100$, and
\item $\varepsilon(\alpha_7)=(1100,1200,2200,2100)$ and 
      $\varepsilon(\alpha_8)=(0100,0200)$. 
\end{itemize}
And in the concatenation
$$
\varepsilon(\alpha_1)\circ\overline{\varepsilon(\alpha_2)}
\circ\varepsilon(\alpha_3)\circ\overline{\varepsilon(\alpha_4)}
\circ\varepsilon(\alpha_5)\circ\overline{\varepsilon(\alpha_6)}
\circ\varepsilon(\alpha_7)\circ\overline{\varepsilon(\alpha_8)}
\circ\varepsilon(\alpha_9)\circ\overline{\varepsilon(\alpha_{10})}
\circ\varepsilon(\alpha_{11})\circ\overline{\varepsilon(\alpha_{12})}
\circ\varepsilon(\alpha_{13}),
$$
in the transition from $\varepsilon(\alpha_7)$ to $\overline{\varepsilon(\alpha_8)}$
we found $2100$ followed by $0200$ which differ in more than one position.

\subsection{The case of $q$ odd}

A way to overcome the previous difficulties is to 
consider the partition of ${\mathcal F_{n}^{(k)}}$ as in its definition (\ref{f_nk}).
For $j=1,2,\ldots,k$, let  $\alpha_i^{(j)}$ be the $i$-th element in 
the list $0^{j-1}1\cdot\overline{\mathcal F_{n-j}^{(k)}}$, and 

$$\Gamma_j=
\varepsilon(\alpha_1^{(j)})\circ\overline{\varepsilon(\alpha_2^{(j)})}
\circ\varepsilon(\alpha_3^{(j)})\circ\overline{\varepsilon(\alpha_4^{(j)})}\circ
\ldots
\circ\varepsilon'(\alpha_{f_{n+k-j}^{(k)}}^{(j)}),
$$
where $\varepsilon'(\alpha_{f_{n+k-j}^{(k)}}^{(j)})=\varepsilon(\alpha_{f_{n+k-j}^{(k)}}^{(j)})$
if $f_{n+k-j}^{(k)}$ is odd and
$\varepsilon'(\alpha_{f_{n+k-j}^{(k)}}^{(j)})=\overline{\varepsilon(\alpha_{f_{n+k-j}^{(k)}}^{(j)})}$
if $f_{n+k-j}^{(k)}$ is even. Let define $\mathcal F_{n,q}^{(k)}$ as
$$\mathcal F_{n,q}^{(k)}
=
\Gamma_1\circ\Gamma_2\circ\ldots\circ\Gamma_k,
$$
and clearly $\mathcal F_{n,q}^{(k)}$ is a list for the set $F_{n,q}^{(k)}$ 
and the next proposition proves that it is a Gray code.

\begin{Proposition}
If $q$ is odd, then the list
$\mathcal F_{n,q}^{(k)}
=
\Gamma_1\circ\Gamma_2\circ\ldots\circ\Gamma_k
$
is a Gray code list with Hamming distance 1.
\end{Proposition}

\noindent
{\it Proof.}
We have to prove the following:
\begin{enumerate}
\item $\Gamma_j$ is a Gray code list, for each $j=1,2,\ldots,k$, with Hamming distance 1;
\item $d_H(\last(\Gamma_j),\first(\Gamma_{j+1}))=1$, for $j=1,2,\ldots,k-1$.
\end{enumerate}
By Proposition \ref{epsilon} it follows that $\varepsilon(\alpha_i^{(j)})$ is a Gray code,
and for the point 1 we have to check
that, if $i$ is odd, $d_H(\last(\varepsilon(\alpha_i^{(j)})),\first(\overline{\varepsilon(\alpha_{i+1}^{(j)})}))=1$
and that, if $i$ is even, $d_H(\last(\overline{\varepsilon(\alpha_i^{(j)})}),\first(\varepsilon(\alpha_{i+1}^{(j)})))=1$.

\begin{itemize}
\item When $i$ is odd we observe that, for some $v$ and $w$,
$$
\begin{array}{ccl}
\alpha_i^{(j)}&=&0^{j-1}1v,\ \mbox{and}\\
\\
\alpha_{i+1}^{(j)}&=&0^{j-1}1w
\end{array}
$$
where $v$ and $w$ differ in a single position 
since
$\alpha_i^{(j)}$ and $\alpha_{i+1}^{(j)}$ are two consecutive 
binary strings in $\mathcal F_n^{(k)}$,
which is a Gray code list.

Let $t=|v|_1$, and since $q-1$ is even, by Proposition \ref{first_last}
it follows that $\last(\mathcal G_{t+1}^{q-1}+1)=(q-1)1^t$ 
which occurs in the last element of the expansion of $\alpha_i^{(j)}$.
Therefore, $\last(\varepsilon(\alpha_i^{(j)}))=0^{j-1}(q-1)v.$

Now, $\first(\overline{\varepsilon(\alpha_{i+1}^{(j)})})$ is equal
to $\last({\varepsilon(\alpha_{i+1}^{(j)})})$, which as previously,
is equal in turn to $0^{j-1}(q-1)w$.

Since $v$ and $w$ differ in a single position, so do 
$\last(\varepsilon(\alpha_i^{(j)}))$
and $\first(\overline{\varepsilon(\alpha_{i+1}^{(j)})})$.

\item If $i$ is even, and since
$\last(\overline{\varepsilon(\alpha_i^{(j)})})=\first(\varepsilon(\alpha_i^{(j)}))$,
by the definition of expansion it follows that
$$ \begin{array}{ccl}
\first(\varepsilon(\alpha_i^{(j)}))&=&\alpha_i^{(j)},\ \mbox{and}\\
\\
\first(\varepsilon(\alpha_{i+1}^{(j)}))&=&\alpha_{i+1}^{(j)}.
\end{array}
$$

Since $\alpha_i^{(j)}$ and $\alpha_{i+1}^{(j)}$ are two consecutive 
strings, their Hamming distance is 1.
\end{itemize}

\noindent For the second point, we have
$$
\begin{array}{rcl}
\first(\Gamma_{j+1})&=&\first(\varepsilon(\alpha_1^{(j+1)})),\ \mbox{and}\\
\\
\last(\Gamma_j)&=&\last(\varepsilon'(\alpha_{f_{n+k-j}^{(k)}}^{(j)})),
\end{array}
$$
where, for some $w'$ and $w''$,
$$
\begin{array}{rcl}
\alpha_1^{(j+1)}&=&0^j1w',\ \mbox{and}\\
\\
\alpha_{f_{n+k-j}^{(k)}}^{(j)}&=&0^{j-1}1w''.
\end{array}
$$
Since $\alpha_{f_{n+k-j}^{(k)}}^{(j)}$ and $\alpha_1^{(j+1)}$ are two consecutive 
strings in $\mathcal F_n^{(k)}$, their Hamming distance is $1$, and thus $w''=1w'$.
Two cases can occur:
\begin{itemize}
\item if $f_{n+k-j}^{(k)}$ is even, then
$
\last(\Gamma_j)=\last(\varepsilon'(\alpha_{f_{n+k-j}^{(k)}}^{(j)}))=\last(\overline{\varepsilon(\alpha_{f_{n+k-j}^{(k)}}^{(j)})})
=\first(\varepsilon(\alpha_{f_{n+k-j}^{(k)}}^{(j)}))=\alpha_{f_{n+k-j}^{(k)}}^{(j)}.
$
Moreover,
$
\first(\Gamma_{j+1})=\first(\varepsilon(\alpha_1^{(j+1)}))=\alpha_1^{(j+1)}.
$
\item if $f_{n+k-j}^{(k)}$ is odd, then
$
\last(\Gamma_j)=\last(\varepsilon(\alpha_{f_{n+k-j}^{(k)}}^{(j)}))=0^{j-1}(q-1)w''$
and
$
\first(\Gamma_{j+1})=\first(\varepsilon(\alpha_1^{(j+1)}))=\alpha_1^{(j+1)}=0^j1w'=0^jw''
$.

\end{itemize}

In any case, $d_H(\last(\Gamma_j),\first(\Gamma_{j+1}))=1$.

Therefore the proof
is concluded and $\mathcal F_{n,q}^{(k)}$ is a Gray code list with Hamming distance 1.
\qed

\medskip

It is easy to see that, generally, when $q$ is even, the construction given
in the previous proposition does not yield a Gray code.

For the sake of clearness, we 
illustrate the previous construction for the Gray code when
$n=4$, $k=3$ and $A=\{0,1,2\}$. We have:
$$
\begin{array}{rcl}
\mathcal F_4^{(3)}&=&\small{(\underbrace{1001,1011,1010,1110,1111,1101,1100}_{\alpha_i^{(1)}\in 1\cdot\overline{\mathcal F_3^{(3)}}},
\underbrace{0100,0101,0111,0110}_{\alpha_i^{(2)}\in 01\cdot\overline{\mathcal F_2^{(3)}}},\underbrace{0010,0011}_
{\alpha_i^{(3)}\in 001\cdot\overline{\mathcal F_1^{(3)}}});}\\

\mathcal G_0^2\oplus 1&=&\lambda;\\
\\
\mathcal G_1^2\oplus 1&=&(1,2);\\
\\
\mathcal G_2^2\oplus 1&=&(11,12,22,21);\\
\\
\mathcal G_3^2\oplus 1&=&(111,112,122,121,221,222,212,211);\\
\\
\mathcal G_4^2\oplus 1&=&(1111,1112,1122,1121,1221,1222,1212,1211,2211,2212,2222,2221,\\
            & &  2121,2122,2112,2111);\\

\end{array}
$$

$$
\begin{array}{rcl}

\Gamma_1&=&(
\underbrace{1001,1002,2002,2001,}_{\varepsilon(\alpha_1^{(1)})}
\underbrace{2011,\ldots,1011,}_{\overline{\varepsilon(\alpha_2^{(1)})}}
\underbrace{1010,\ldots,2010,}_{\varepsilon(\alpha_3^{(1)})}
\underbrace{2110,\ldots,1110,}_{\overline{\varepsilon(\alpha_4^{(1)})}}\\
& &
\underbrace{1111,\ldots,2111,}_{\varepsilon(\alpha_5^{(1)})}
\underbrace{2101,\ldots,1101,}_{\overline{\varepsilon(\alpha_6^{(1)})}}
\underbrace{1100,\ldots,2100}_{\varepsilon(\alpha_7^{(1)})});\\
\\
\Gamma_2&=&(
\underbrace{0100,0200,}_{\varepsilon(\alpha_1^{(2)})}
\underbrace{0201,\ldots,0101,}_{\overline{\varepsilon(\alpha_2^{(2)})}}
\underbrace{0111,\ldots,0211,}_{\varepsilon(\alpha_3^{(2)})}
\underbrace{0210,\ldots,0110}_{\overline{\varepsilon(\alpha_4^{(2)})}});\\
\\
\Gamma_3&=&(
\underbrace{0010,0020,}_{\varepsilon(\alpha_1^{(3)})}
\underbrace{0021,\ldots,0011}_{\overline{\varepsilon(\alpha_2^{(3)})}}).\\

\end{array}
$$

The reader can easily check that $\mathcal F_{4,3}^{(3)}=\Gamma_1\circ\Gamma_2\circ\Gamma_3$
is a Gray code with Hamming distance~1.

\subsection{The case of $q$ even}
In this case the construction of
a Gray code list for $F_{n,q}^{(k)}$ is straightforward, and 
based on the discussion at the beginning of this section: just consider the
expansions of the binary strings  $\alpha_i$ in $\mathcal F_n^{(k)}$, 
for $i=1,2,\ldots,f_{n+k}^{(k)}$, and concatenate them, taking $\varepsilon(\alpha_i)$
and $\overline{\varepsilon(\alpha_{i+1})}$ alternatively, 
as in expression (\ref{Gent_Alter}).
The next proposition shows that the obtained list is a Gray code.
 
\begin{Proposition}
If $q$ is even, then the list
$$
\mathcal{F}_{n,q}^{(k)}=\varepsilon(\alpha_1)\circ\overline{\varepsilon(\alpha_2)}
\circ\ldots\circ\varepsilon'(\alpha_{f_{n+k}^{(k)}}),
$$
where $\varepsilon'(\alpha_{f_{n+k}^{(k)}})=\varepsilon(\alpha_{f_{n+k}^{(k)}})$
if $f_{n+k}^{(k)}$ is odd and
$\varepsilon'(\alpha_{f_{n+k}^{(k)}})=\overline{\varepsilon(\alpha_{f_{n+k}^{(k)}})}$
if $f_{n+k}^{(k)}$ is even,
is a Gray code list with Hamming distance 1.
\end{Proposition}
\noindent
{\it Proof.}
By Proposition \ref{epsilon}, each $\varepsilon(\alpha_i)$ is a Gray code list,
and, we have to check that, if $i$ is odd,
$d_H(\last(\varepsilon(\alpha_i)),\first(\overline{\varepsilon(\alpha_{i+1}}))=1$
 and that, if $i$ is even,
 $d_H(\last(\overline{\varepsilon(\alpha_i)}),\first(\varepsilon(\alpha_{i+1})))=1$.

\begin{itemize}
\item In the first case ($i$ is odd), by Proposition \ref{first_last}
and the definition of expansion, it follows that
in $\last(\varepsilon(\alpha_i))$ and in $\last(\varepsilon(\alpha_{i+1}))$
the symbols different from $0$ are
equal to $q-1$. For example, if $\alpha_i=10110$, then
$\last(\varepsilon(\alpha_i))=(q-1)0(q-1)(q-1)0$.
Moreover, since $\alpha_i$ and $\alpha_{i+1}$ are two consecutive strings of
$\mathcal F_n^{(k)}$, $d_H(\alpha_i,\alpha_{i+1})=1$, and so
$d_H(\last(\varepsilon(\alpha_i)), \last(\varepsilon(\alpha_{i+1})))=1$.

\item If $i$ is even, we observe that
$$
\begin{array}{rclc}
\last(\overline{\varepsilon(\alpha_i)})&=&\first(\varepsilon(\alpha_i))&=\alpha_i,\ \mbox{and}\\
\\
\first(\varepsilon(\alpha_{i+1}))&=&\alpha_{i+1}.
\end{array}
$$
Since $d_H(\alpha_i,\alpha_{i+1})=1$, then
$d_H(\last(\overline{\varepsilon(\alpha_i)}),\first(\varepsilon(\alpha_{i+1})))=1$.
\end{itemize}
\qed

\section{Conclusions and further developments}
In this paper we propose a trace partitioned Gray code for the $q$-ary $k$-generalized 
Fibonacci strings of length $n$, where
the Hamming distance between two successive strings is $1$. 
Our constructions are based on the expansion of an existing Gray code 
when $q=2$.
A consequence of the expansion technique is that 
our  Gray code has the following property:
if we replace each non-zero symbol in each string by $1$,
and `collapse' the obtained list by  
keeping one copy of each binary strings, then
the existing Gray code for $q=2$ is obtained.

The investigation on the existence of Gray codes
for strings on a $q$-ary alphabet avoiding a general consecutive
pattern has already been studied: in \cite{Sq} the author gives
such a Gray code only when $q$
even, while the case of $q$ odd is left open. 
Our Gray code deals with the avoidance of a particular pattern
but works for any $q$, and 
an interesting development could be a deeper
investigation to check if this constructions can be applied to a
general consecutive pattern in order to solve the open question in \cite{Sq}.
Also, it would be of interest to implement our Gray code definition
into an efficient generating algorithm for the set of underlying strings.


\begin{thebibliography}{99}

\bibitem{E} Er, M. C.:
On generating the $N$-ary reflected Gray code.
IEEE Trans. Comput. {\bfseries33} (1984) 739--741.

\bibitem{G} Gray, F.:
Pulse Code Communication
U.S. Patent 2 632 058 (March 17, 1953).

\bibitem{H} Hamming, R. W.:
Error detecting and error correcting codes.
Bell System Tech. J. {\bfseries29} (1950) 147--160.

\bibitem{K} Knuth, D. E.:
The Art of Computer Programming.
Vol. 3, ``Sorting and Searching", Addison-Wesley, Reading, MA, 1973.

\bibitem{L} Levesque, C.:
On $m$-th order linear recurrences.
Fibonacci Quart. {\bfseries 23} (1985) 290--293.

\bibitem{Rus} 
Ruskey, F.: 
Simple combinatorial Gray codes constructed by reversing sublists, 
4th ISAAC, 
Lecture Notes in Computer Science,  {\bfseries 762} (1993) 201-208. 

\bibitem{Sh} Shork, M.: The $r$-generalized Fibonacci numbers and Polynomial coefficients.
Int. J. Contemp. Math. Sciences. {\bfseries3} (2008) 1157--1163.

\bibitem{Sq} Squire, M. B.: Gray codes for $A$-free strings.
Electron. J. Combin. {\bfseries3} (1996) \#R17.

\bibitem{V1} Vajnovszki, V.:
A loopless generation of bitstrings without $p$ consecutive ones.
Discrete Math. Theor. Comput. Sci. Springer 2001, 227--240.

\bibitem{Wi} Williamson, S. G.:
Combinatorics for computer science.
Computer Science Press, Rockville, Maryland, 1985.

\end{thebibliography}
\end{document}